\documentclass[a4paper,11pt]{amsart}
\usepackage{amsmath,amsthm,amssymb,enumerate}
\usepackage{graphicx}

\numberwithin{equation}{section}
\numberwithin{figure}{section}

\newtheorem{thm}{Theorem}[section]

\newtheorem{definition}[thm]{{Definition}}
\newtheorem{remark}[thm]{{Remark}}
\newtheorem{example}[thm]{\rm{Example}}

\newenvironment{df}{\begin{definition}\rm}{\end{definition}}
\newenvironment{rem}{\begin{remark}\rm}{\end{remark}}

\title{Surface links with free abelian groups}
\author{Inasa Nakamura}
\address{
Department of Mathematics, 
Gakushuin University, 
1-5-1 Mejiro Toshima-ku Tokyo, 171-8588 JAPAN} 

\email{inasa@math.gakushuin.ac.jp}



\begin{document}
\maketitle

\begin{abstract}
It is known that if a classical link group is a free abelian group, then its rank is at most two. It is also known that a $k$-component 2-link group ($k>1$) is not free abelian. 
In this paper, we give examples of $T^2$-links each of whose link groups is a free abelian group of rank three or four. Concerning the $T^2$-links of rank three, we determine the triple point numbers and we see that their link types are infinitely many. 
\end{abstract}

\section{Introduction}
A {\it classical link} is the image of a smooth embedding of a disjoint union of circles into the Euclidean 3-space $\mathbb{R}^3$. 
The {\it link group} is the fundamental group of the link exterior. 
It is known \cite[Theorem 6.3.1]{Kawauchi} that if a classical link group is a free abelian group, then its rank is at most two. 
A {\it surface link} is the image of a smooth embedding of a closed surface into the Euclidean 4-space $\mathbb{R}^4$. 
A {\it $2$-link} (resp. {\it $T^2$-link}) is a surface link whose components are homeomorphic to 2-spheres (resp. tori). 
It is known \cite[Chapter 3, Corollary 2]{Hillman} that a $k$-component 2-link group for $k>1$ is not a free abelian group.  
The aim of this paper is to give concrete examples of $T^2$-links whose link groups are free abelian. 

It is known (see Remark \ref{remark1}) that a $T^2$-link called a \lq\lq Hopf 2-link" \cite{CKSS01} has a free abelian group of rank two. We give $T^2$-links with a free abelian group of rank three (Theorem \ref{Thm1}). 
We also give a $T^2$-link with a free abelian group of rank four (Theorem \ref{Thm2}). 
These $T^2$-links are \lq\lq torus-covering $T^2$-links", which are $T^2$-links in the form of unbranched coverings over the standard torus. 

Further we study the $T^2$-links given in Theorem \ref{Thm1} i.e. $T^2$-links each of whose link groups is a free abelian group of rank three. 
We determine the triple point number of each $T^2$-link (Theorem \ref{Thm3}), by which we can see that their link types are infinitely many. 
The triple point number of each $T^2$-link is a multiple of four, and it is realized by a surface diagram in the form of a covering over the torus. 
For other examples of surface links (not necessarily orientable) which realize large triple point numbers, see \cite{COS, Kamada92-2, Kamada-Oshiro, N2, Oshiro10, Satoh}. 
 
The paper is organized as follows. 
In Section \ref{section1}, we review the definition of a torus-covering $T^2$-link, and we review a formula how to calculate its link group. 
In Section \ref{section2}, we show Theorems \ref{Thm1} and \ref{Thm2}. In Section \ref{TriplePoint}, we show Theorem \ref{Thm3}. 

\section{A torus-covering $T^2$-link and its link group} \label{section1}
 
In this section, we give the definition of a torus-covering $T^2$-link $\mathcal{S}_m(a,b)$, which is determined from a pair of commuting $m$-braids $a$ and $b$ called basis braids. 
For the definition of a torus-covering link whose component might be of genus more than one, see \cite{N}. 
We can compute the link group of $\mathcal{S}_m(a,b)$ by using Artin's automorphism associated with $a$ or $b$ \cite{N}. 

\subsection{} \label{0913-1}
Let $T$ be the standard torus in $\mathbb{R}^4$, i.e. the boundary of an unknotted solid torus in $\mathbb{R}^3 \times \{0\} \subset \mathbb{R}^4$. 
Let $N(T)$ be a tubular neighborhood of $T$ in $\mathbb{R}^4$. 

\begin{df} \label{Def2-1} 
A {\it torus-covering $T^2$-link} is a surface link $F$ in $\mathbb{R}^4$ such that $F$ is embedded in $N(T)$ and $p |_{F} \,:\, F \to T$ is an unbranched covering map, where $p \,:\, N(T) \to T$ is the natural projection. 
\end{df}
 
Let us consider a torus-covering $T^2$-link $F$. 
Let us fix a point $x_0$ of $T$, and take a meridian $\mathbf{m}$ and a longitude $\mathbf{l}$ of $T$ with the base point $x_0$. A {\it meridian} is an oriented simple closed curve on $T$ which bounds a 2-disk in the solid torus whose boundary is $T$ and which is not null-homologous in $T$. A {\it longitude} is an oriented simple closed curve on $T$ which is null-homologous in the complement of the solid torus in the three space $\mathbb{R}^3 \times \{0\}$ and which is not null-homologous in $T$. The intersections $F \cap p^{-1}(\mathbf{m})$ and $F \cap p^{-1}(\mathbf{l})$ are closures of classical braids. 
Cutting open the solid tori at the 2-disk $p^{-1}(x_0)$, we obtain a pair of classical braids. We call them {\it basis braids} \cite{N}. 
The basis braids of a torus-covering $T^2$-link are commutative, and for any commutative braids $a$ and $b$, there exists a unique torus-covering $T^2$-link with basis braids $a$ and $b$ \cite[Lemma 2.8]{N}. 
For commutative $m$-braids $a$ and $b$, we denote by $\mathcal{S}_m(a,b)$ the torus-covering $T^2$-link with basis $m$-braids $a$ and $b$. 
 
\subsection{} \label{0913-2}
We can compute the link group of a torus-covering $T^2$-link $\mathcal{S}_m(a,b)$ \cite{N}. 
As preliminaries, we will give the definition of Artin's automorphism (see \cite{Kamada3}). Let $c$ be an $m$-braid in a cylinder $D^2 \times [0,1]$, and let $Q_m$ be the starting point set of $c$. Let $\{h_u\}_{u \in [0,1]}$ be an isotopy of $D^2$ rel $\partial D^2$ such that $\cup_{u \in [0,1]} h_u (Q_m) \times \{u\}=c$. 
 Let $\mathcal{A}^c \,:\, (D^2, Q_m) \rightarrow (D^2, Q_m)$ be the terminal map $h_1$, and consider the induced map $\mathcal{A}^c_* \,:\, \pi_1(D^2-Q_m) \rightarrow \pi_1(D^2-Q_m)$. It is known \cite{Artin} that $\mathcal{A}^c_*$ is uniquely determined from $c$. We call $\mathcal{A}^c_*$ {\it Artin's automorphism} associated with $c$. Note that $\pi_1(D^2-Q_m)$ is naturally isomorphic to the free group $F_m$ generated by the standard generators $x_1, x_2, \ldots, x_m$ of $\pi_1(D^2-Q_m)$. 
By $\mathcal{A}^c_*$, the braid group $B_m$ acts on $\pi_1(D^2-Q_m)$. It is presented by 
\begin{gather*}
\mathcal{A}^{\sigma_i}_*(x_j)=\begin{cases}
                        x_j x_{j+1} x_j^{-1} & \text{if $j = i$} \\
                        x_{j-1} & \text{if $j = i+1$} \\
                        x_j & \text{otherwise} 
\end{cases}\\
\tag*{\text{and}}
\mathcal{A}^{\sigma_i^{-1}}_*(x_j)=\begin{cases}
                        x_{j+1} & \text{if $j = i$} \\
                        x_{j}^{-1} x_{j-1} x_{j} & \text{if $j = i+1$} \\
                        x_j & \text{otherwise}
\end{cases}
\end{gather*}
where $i=1,2,\ldots,m-1$ and $j=1,2,\ldots,m$. 

It is known \cite[Proposition 3.1]{N} that the link group of $\mathcal{S}_m(a,b)$ is presented by 
\begin{equation*} 
\pi_1(\mathbb{R}^4-\mathcal{S}_m(a,b))=\langle \, x_1 \,, \ldots,\, x_m \mid 
x_j=\mathcal{A}^a_*(x_j)=\mathcal{A}^b_* (x_j), \text{ for $ j =1,2,\ldots,m$} \, \rangle. 
\end{equation*}

\section{$T^2$-links whose link groups are free abelian} \label{section2}
In this section we show Theorems \ref{Thm1} and \ref{Thm2}: There are torus-covering $T^2$-links with a free abelian group of rank three (Theorem \ref{Thm1}) or four (Theorem \ref{Thm2}). 

\begin{rem} \label{remark1}
A {\it Hopf 2-link} \cite{CKSS01} is a $T^2$-link which is the product of a classical Hopf link in $B^3$ with $S^1$, embedded into $\mathbb{R}^4$ via an embedding of $B^3 \times S^1$ into $\mathbb{R}^4$, where $B^3$ is a 3-ball and $S^1$ is a circle. 
There are two link types according to the embedding of $B^3 \times S^1$, called a standard Hopf 2-link and a twisted Hopf 2-link \cite{CKSS01}. A standard (resp. twisted) Hopf 2-link is the spun $T^2$-link (resp. the turned spun $T^2$-link) of a classical Hopf link \cite{Livingston, Boyle88, Boyle}. 
It is known \cite{Livingston, Boyle88, Boyle} that the link group of the spun $T^2$-link or the turned spun $T^2$-link of a classical link $L$ is isomorphic to the classical link group of $L$. Thus we can see that a Hopf 2-link has a free abelian link group of rank two. 
\end{rem}
 
Let $\sigma_1, \sigma_2, \ldots, \sigma_{m-1}$ be the standard generators of $B_m$. 

\begin{thm} \label{Thm1}
The link group of $\mathcal{S}_3(\sigma_1^2 \sigma_2^{2n}, \Delta)$ is a free abelian group of rank three, where $n$ is an integer and $\Delta$ is a full twist of a bundle of three parallel strings. 
\end{thm}

\begin{proof}
Put $S_n=\mathcal{S}_3(\sigma_1^2 \sigma_2^{2n}, \Delta)$. 
Let us compute the link group $G_n=\pi_1(\mathbb{R}^4-S_n)$ 
by applying \cite[Proposition 3.1]{N}. 
Let $x_1$, $x_2$ and $x_3$ be the generators. 
Then the relations concerning the basis braid $\sigma_1^2 \sigma_2^{2n}$ are
\begin{eqnarray}
x_1 x_2 &=& x_2 x_1, \label{4-1} \\
(x_2 x_3)^{|n|} &=& (x_3 x_2)^{|n|}.  
\end{eqnarray}
The other relations concerning the other basis braid $\Delta$ are 
\begin{eqnarray*}
x_1 &=& (x_1 x_2 x_3 ) x_1 (x_1 x_2 x_3)^{-1}, \\
x_2 &=& (x_1 x_2 x_3 ) x_2 (x_1 x_2 x_3)^{-1}, \\
x_3 &=& (x_1 x_2 x_3 ) x_3 (x_1 x_2 x_3)^{-1}, 
\end{eqnarray*}
which are 
\begin{eqnarray}
   x_1 x_2 x_3  &=& x_2 x_3 x_1, \label{4-3} \\
   x_2 (x_1 x_2 x_3 ) &=& (x_1 x_2 x_3) x_2, \label{4-4} \\
 x_3 x_1 x_2  &=& x_1 x_2 x_3 .  
\end{eqnarray}
By (\ref{4-1}), (\ref{4-3}) is deformed to $x_2 x_1 x_3=x_2 x_3 x_1$; hence  
\begin{equation} \label{4-6}
  x_1 x_3=x_3 x_1. 
\end{equation}
Similarly, by (\ref{4-4}) and (\ref{4-1}), 
\begin{equation} \label{4-6-2}
 x_2 x_3=x_3 x_2. 
 \end{equation}
We can see that all the relations are generated by 
the three relations (\ref{4-1}), (\ref{4-6}) and (\ref{4-6-2}). 
 Thus we have 
 \begin{eqnarray*}
 G_n &=& \langle \,x_1,\, x_2,\, x_3 \mid x_1 x_2=x_2 x_1, \, 
x_2 x_3=x_3 x_2, \, x_3 x_1 =  x_1 x_3 \, \rangle, 
 \end{eqnarray*}
which is a free abelian group of rank three. 
\end{proof}

\begin{thm} \label{Thm2}
The link group of $\mathcal{S}_4(\sigma_1^2 \sigma_2^2 \sigma_3^2, \Delta)$ is a free abelian group of rank four, where $\Delta$ is a full twist of a bundle of $4$ parallel strings. 
\end{thm}
\begin{proof}
 Similarly to the proof of Theorem \ref{Thm1}, by \cite[Proposition 3.1]{N}, for generators $x_1$, $x_2$, $x_3$ and $x_4$, we have the following relations:
 \begin{equation}\label{i}
 x_i x_{i+1} = x_{i+1} x_i, 
 \end{equation}
  where $i=1$,$2$,$3$, and 
  \begin{equation}\label{eq}
x_i = (x_1 x_2 x_3 x_4) x_i (x_1 x_2 x_3 x_4)^{-1},
\end{equation}
 where $i=1$,$2$, $3$, $4$.  
Using $x_1 x_2=x_2 x_1$ and $x_3 x_4=x_4 x_3$ of (\ref{i}), the latter four relations (\ref{eq}) are deformed as follows: 
\begin{eqnarray}
 x_1 x_3 x_4 &=& x_3 x_4 x_1, \\
  x_2 x_3 x_4 &=& x_3 x_4 x_2, \label{4-7} \\
 x_3 x_1 x_2 &=& x_1 x_2 x_3, \label{4-8} \\
 x_4 x_1 x_2 &=& x_1 x_2 x_4. \label{4-9}
\end{eqnarray}
By $x_2 x_3=x_3 x_2$ of (\ref{i}), (\ref{4-7}) is deformed to 
$x_3 x_2 x_4=x_3 x_4 x_2$; hence  
\begin{equation} \label{5-1}
  x_2 x_4=x_4 x_2. 
\end{equation}
Similarly, by $x_2 x_3=x_3 x_2$ of (\ref{i}) and (\ref{4-8}), 
\begin{equation} \label{5-2}
x_3 x_1=x_1 x_3, 
\end{equation}
and by (\ref{5-1}) and (\ref{4-9}), 
\begin{equation} \label{5-3}
 x_4 x_1=x_1 x_4. 
 \end{equation}
We can see that all the relations are generated by 
the relations (\ref{i}), (\ref{5-1}), (\ref{5-2}) and (\ref{5-3}). 
Thus the link group is a free abelian group of rank four. 
\end{proof}

\section{The triple point numbers of the $T^2$-links 
with a free abelian group of rank three} \label{TriplePoint}
The {\it triple point number} of a surface link $F$ is the minimal number of triple points among all the surface diagrams of $F$. 
 In this section we study the $T^2$-links given in Theorem \ref{Thm1} i.e. $T^2$-links each of whose link group is a free abelian group of rank three. 

\begin{thm}\label{Thm3}
The triple point number of $S_n=\mathcal{S}_3(\sigma_1^2 \sigma_2^{2n}, \Delta)$ given in Theorem \ref{Thm1} is $4n$ for $n>0$ and $4(1-n)$ for $n \leq 0$. Further it is realized by a surface diagram in the form of a covering over $T$, in other words, by a $3$-chart on $T$ which presents $S_n$. Thus $T^2$-links with a free abelian group of rank three are infinitely many. 
\end{thm}
\noindent
Here, a $3$-chart \cite{Kamada3} is a finite graph with certain additional data, which we review in Section \ref{1130-1}. 

This section is organized as follows. 
In Section \ref{1130-1}, we review a surface diagram and an $m$-chart on $T$ which presents a torus-covering $T^2$-link (see \cite{N, Kamada3}). 
In Section \ref{1126-1}, we review the result of \cite{N2} which gives lower bounds of triple point numbers. 
In Section \ref{1207-1}, we prove Theorem \ref{Thm3}.
 
\subsection{Surface diagrams and $m$-charts presenting  torus-covering $T^2$-links} \label{1130-1}
The notion of an $m$-chart on a 2-disk was introduced by Kamada \cite{Kamada92} (see also \cite{Kamada3}) to present a surface braid i.e. a 2-dimensional braid in a bi-disk (see \cite{Rudolph, Kamada3}). 
An $m$-chart on a disk is obtained from the singularity set of a surface diagram of a surface braid.
By a minor modification, we can define an $m$-chart on $T$ presenting a torus-covering link \cite{N}.

For a torus-covering $T^2$-link $F$, we consider a surface diagram in the form of a covering over the torus, as in Section \ref{sec3-1-1}. 
Given $F$, we obtain such a surface diagram $D$, and from $D$ we obtain a graph called an $m$-chart on $T$ (without black vertices).  Conversely, an $m$-chart on $T$ without black vertices presents such a surface diagram and hence a torus-covering $T^2$-link. 

\subsubsection{Surface diagrams}\label{sec3-1-1}

\begin{figure}
\begin{center}
\includegraphics*{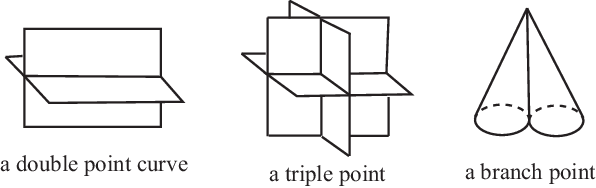}
\end{center}
\caption{The singularity of a surface diagram.}
\label{0215-1}
\end{figure}

We review a surface diagram of a surface link $F$ (see \cite{CKS}). For a projection $\pi \,:\, \mathbb{R}^4 \to \mathbb{R}^3$, the closure of the self-intersection set of $\pi(F)$ is called the singularity set. Let $\pi$ be a generic projection, i.e. the singularity set of the image $\pi(F)$ consists of double points, isolated triple points, and isolated branch points; see Fig. \ref{0215-1}. The closure of the singularity set forms a union of immersed arcs and loops, which we call double point curves. Triple points (resp. branch points) form the intersection points (resp. the end points) of the double point curves. A {\it surface diagram} of $F$ is the image $\pi(F)$ equipped with over/under information along each double point curve with respect to the projection direction. 

Throughout this paper, we consider the surface diagram of a torus-covering $T^2$-link $F$ by the projection which projects $N(T)=I \times I \times T$ to $I \times T$ for an interval $I$, where we identify $N(T)$ with $I \times I \times T$ in such a way as follows. Since $T$ is the boundary of the standard solid torus in $\mathbb{R}^3 \times \{0\}$, the normal bundle of $T$ in $\mathbb{R}^3 \times \{0\}$ is a trivial bundle. We identify it with $I \times T$. Then we identify $N(T)$ with $I \times I \times T$, where the second $I$ is an interval in the fourth axis of $\mathbb{R}^4$. Perturbing $F$ if necessary, we can assume that this projection is generic. We call this surface diagram {\it the surface diagram of $F$ in the form of a covering over the torus}. 

\subsubsection{From surface diagrams to $m$-charts on $T$}
Given a torus-covering $T^2$-link $F$, 
we obtain a graph on $T$ from the surface diagram in the form of a covering over the torus, as follows. 
Now we have $\mathrm{Sing}(\pi(F))$ in $I \times T$. By the definition of a torus-covering $T^2$-link, $\mathrm{Sing}(\pi(F))$ consists of  double point curves and triple points, and no branch points.
We can assume that the singular set of the image of $\mathrm{Sing}(\pi(F))$ by the projection to $T$ consists of a finite number of double points such that the preimages belong to double point curves of $\mathrm{Sing}(\pi(F))$. 
Thus the image of $\mathrm{Sing}(\pi(F))$ by the projection to $T$ forms a finite graph $\Gamma$ on $T$ such that the degree of its vertex is either $4$ or $6$. An edge of $\Gamma$ corresponds to a double point curve, and a vertex of degree $6$ corresponds to a triple point.  

For such a graph $\Gamma$ obtained from the surface diagram, we give orientations and labels to the edges of $\Gamma$, as follows. Let us consider a path $l$ in $T$ such that $l \cap \Gamma$ is a point $P$ of an edge $e$ of $\Gamma$. Then $F \cap p^{-1} (l)$ is a classical $m$-braid with one crossing in $p^{-1}(l)$ such that $P$ corresponds to the crossing of the $m$-braid. Let $\sigma_{i}^{\epsilon}$ ($i \in \{1,2,\ldots, m-1\}$, $\epsilon \in \{+1, -1\}$) be the presentation of $F \cap p^{-1}(l)$. Then label the edge $e$ by $i$, and moreover give $e$ an orientation such that the normal vector of $l$ corresponds (resp. does not correspond) to the orientation of $e$ if $\epsilon=+1$ (resp. $-1$). We call such an oriented and labeled graph an {\it $m$-chart of $F$} (without black vertices). 
\\

 In general, we define an $m$-chart on $T$ as follows. 

\begin{figure}
 \includegraphics*{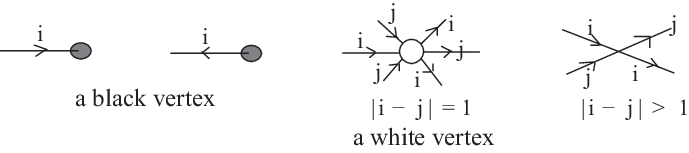}
 \caption{Vertices in an $m$-chart.}
 \label{Fig1-1}
 \end{figure}

\begin{df}
Let $m$ be a positive integer, and let 
$\Gamma$ be a finite graph on $T$.
Then $\Gamma$ is called an {\it $m$-chart on $T$} if 
it satisfies the following conditions: 

\begin{enumerate}[(i)]
 \item Every edge is oriented and labeled by an element of 
       $\{1,2, \ldots, m-1\}$. 
 \item Every vertex has degree $1$, $4$, or $6$.
 \item  The adjacent edges around each vertex are oriented and labeled as shown in 
Fig. \ref{Fig1-1}, 
where we depict a vertex of degree $1$ (resp. 
$6$ by a black vertex (resp. white vertex). 
 \end{enumerate}
 \end{df}
 
 A black vertex presents a branch point; see \cite{Kamada3}.  
When an $m$-chart on $T$ without black vertices is given, we can reconstruct a torus-covering $T^2$-link \cite{N} (see also \cite{Kamada3}). 
\\

\begin{figure}
 \includegraphics*{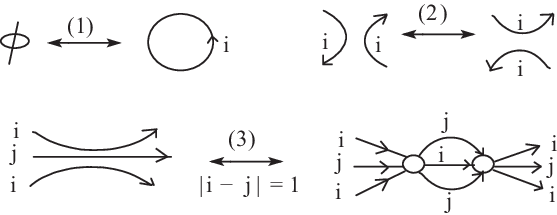}
\caption{CI-moves. We give only several examples.} 
\label{cmove}
 \end{figure}

Two $m$-charts on $T$ are {\it C-move equivalent} \cite{N} (see also \cite{Kamada92, Kamada2, Kamada3}) if they are related by a finite sequence of ambient isotopies of $T$ and CI, CII, CIII-moves. 
We show several examples of CI-moves in Fig. \ref{cmove}; see \cite{Kamada3} for the complete set of CI-moves and CII, CIII-moves. 
For two $m$-charts on $T$, their presenting torus-covering links are equivalent if the $m$-charts are C-move equivalent \cite{N} (see also \cite{Kamada92, Kamada2, Kamada3}). 
 
%

\subsection{Triple point numbers} \label{1126-1}
For a surface link $F$, we denote by $t(F)$ the triple point number of $F$.
It is shown \cite{N2} that for a pure $m$-braid $b$ ($m \geq 3$) and an integer $n$, a lower bound of $t(\mathcal{S}_m(b, \Delta^n))$ is given by using the linking numbers of $\hat{b}$, and for a particular $b$, we can determine the triple point number. Here $\hat{b}$ denotes the closure of $b$.
 
For a pure $3$-braid $b$, it follows from \cite{N2} that we can give a lower bound of $t(\mathcal{S}_3(b, \Delta))$ as follows.  
We define the $i$th component of $\hat{b}$ by the component constructed by the $i$th string of $\hat{b}$ ($i=1,2,3$). For positive integers $i$ and $j$ with $i \neq j$, the {\it linking number} of the $i$th and $j$th components of a classical link $L$, denoted by $\mathrm{Lk}_{i,j}(L)$, is the total number of positive crossings minus the total number of negative crossings of a diagram of $L$ such that the under-arc (resp. over-arc) is from the $i$th (resp. $j$th) component. 
Put $\mu=\sum_{i<j}|\mathrm{Lk}_{i,j}(\hat{b})|$, and put $\nu=\nu_{1,2,3}+\nu_{2,3,1}+\nu_{3,1,2}$, where $\nu_{i,j,k}=\min_{i,j,k} \{ |\mathrm{Lk}_{i,j}(\hat{b})|, |\mathrm{Lk}_{j,k}(\hat{b})| \}$ if $\mathrm{Lk}_{i,j}(\hat{b}) \mathrm{Lk}_{j,k}(\hat{b})>0$ and otherwise zero. 
Then, by \cite{N2}, 
\[
t(\mathcal{S}_3(b, \Delta)) \geq 4(\mu-\nu). 
\]
In particular, let $b$ be a $3$-braid presented by a braid word which is an element of a monoid generated by $\sigma_1^{2}$ and $\sigma_2^{-2}$; note that $b$ is a pure braid.  
Then 
\[
t(\mathcal{S}_3(b, \Delta))=4 \mu, 
\]
and the triple point number is realized by a surface diagram in the form of a covering over the torus \cite{N2}. 
  
\subsection{Proof of Theorem \ref{Thm3}} \label{1207-1}

Put $b=\sigma_1^2 \sigma_2^{2n}$. 
We use the notations given in Section \ref{1126-1}. 
Since $\mathrm{Lk}_{i,j}(\hat{b})=1$ (resp. $n$) if $\{ i,j \}=\{1,2\}$ (resp. $\{2,3\}$) and otherwise zero, we can see that $\mu=1+|n|$. 

Let us consider the case for $n \leq 0$. Since $b$ has the presentation which is an element of a monoid generated by $\sigma_1^2$ and $\sigma_2^{-2}$, $t(S_n)=4 \mu$ by \cite{N2}; thus $t(S_n)=4(1-n)$ ($n \leq 0$), and the triple point number is realized by a surface diagram in the form of a covering over the torus by \cite{N2}. 

Let us consider the case for $n>0$. 
Since $\mathrm{Lk}_{i,j}(\hat{b})=1$ (resp. $n$) if $\{i,j\}=\{1,2\}$ (resp. $\{2,3\}$) and otherwise zero, we can see that $\nu_{i,j,k}=1$ if $(i,j,k)=(1,2,3)$ and zero if $(i,j,k)=(2,3,1)$ or $(3,1,2)$; thus $\nu=1$, and hence $t(S_n) \geq 4(\mu-\nu)=4n$ by \cite{N2}. 

It remains to show that there is a surface diagram of $S_n$ ($n>0$) with $4n$ triple points. 
It suffices to draw a $3$-chart $\Gamma$ on $T$ which presents $S_n$ such that $\Gamma$ has exactly $4n$ white vertices. 
We draw $\Gamma$ which presents $S_n$, and deform it to a $3$-chart with $4n$ white vertices by C-moves, as follows. 
First we draw $\Gamma$ as a $3$-chart which consists of $2n+2$ parts as follows, where we assume that a full twist $\Delta$ has the presentation $\Delta=(\sigma_1\sigma_2\sigma_1)^2$. 
\begin{enumerate}[(i)]
\item
The part of $\Gamma$ with basis braids $\sigma_1$ and 
  $\Delta$. We have two copies. 
 
\item
The part of $\Gamma$ with basis braids $\sigma_2$ and 
  $\Delta$. We have $2n$ copies. 
\end{enumerate}
\begin{figure}
\includegraphics*{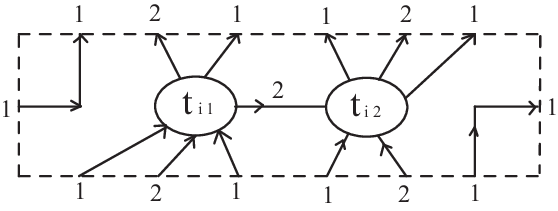}
\caption{White vertices $t_{i1}$ and $t_{i2}$ ($i=1$, $2$).}
\label{Fig3-3}
\end{figure}
\begin{figure}
\includegraphics*{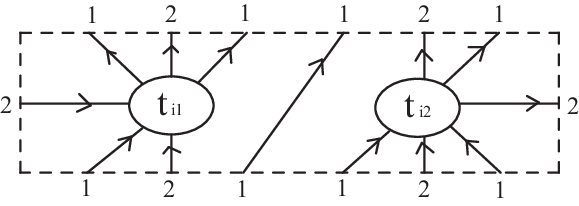}
\caption{White vertices $t_{i1}$ and $t_{i2}$ ($i=3,4,\ldots,2n+2$), for $n>0$.}
\label{Fig3-5}
\end{figure}
We draw the part (i) as in Fig. \ref{Fig3-3} and we denote the white vertices by $t_{i1}$ and $t_{i2}$ as in Fig. \ref{Fig3-3} for the $i$th copy ($i=1,2$). 
We draw the part (ii) as in Fig. \ref{Fig3-5} and we denote the white vertices by $t_{i1}$ and $t_{i2}$ as in Fig. \ref{Fig3-5} for the $(i-2)$th copy ($i=3,4,\ldots,2n+2$). 
There are $4n+4$ white vertices in $\Gamma$. 
Let us apply a CI-move as in Fig. \ref{cmove} (3) to the pair $\{t_{21}, t_{31}\}$ of white vertices in $\Gamma$, and then to the pair $\{t_{(2n+2)2}, t_{12}\}$; then we can eliminate the four white vertices, and the resulting $3$-chart has $4n$ white vertices. Hence $t(S_n)=4n$ ($n>0$), and the triple point number is realized by this $3$-chart on $T$. 
 \qed
 
\begin{figure}
 \includegraphics*{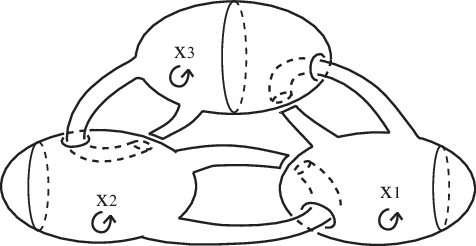}
 \caption{A ribbon $T^2$-link with a free abelian group of rank three.}
 \label{11-0218-1}
 \end{figure}

 \begin{rem}
There is an oriented $T^2$-link as in Fig. \ref{11-0218-1} with a free abelian group of rank three and with the triple point number zero. 
It is a ribbon $T^2$-link (see \cite{CKS} for the definition of a ribbon surface link). 
We briefly show that the link group is free abelian, as follows. 
In the surface diagram, there are six broken sheets (see \cite{CKS}), consisting of three pairs of a sheet attached with $x_i$ and a small disk $D_i$ such that each pair forms the $i$th component of the $T^2$-link ($i=1,2,3$). 
Let us attach $y_i$ to each $D_i$. 
The link group has the presentation with generators $x_i$ and $y_i$ ($i=1,2,3$) and the relations which are given around each double point curve (see \cite{CKS, Kamada3}). 
The singularity set consists of double point curves which form six circles. 
Around each circle in the $i$th component which does not bound $D_i$ ($i=1,2,3$), there are three broken sheets such that one is an over-sheet with $x_i$ and the other two are under-sheets with the same generator $x_{i+1}$, where $x_4=x_1$; together with the orientation, the relation is $x_i= x_{i+1} x_i x_{i+1}^{-1}$, see  \cite{CKS, Kamada3}. 
Around each circle $\partial D_i$ ($i=1,2,3$), there are three broken sheets such that one is an over-sheet with $x_{i+1}$ and the other two are under-sheets with $x_i$ and $y_i$ respectively, where $x_4=x_1$; together with the orientation, the relation is $y_i= x_{i+1} x_i x_{i+1}^{-1}$, see  \cite{CKS, Kamada3}. 
Thus the link group is a free abelian group of rank three. 
\end{rem}

\section*{Acknowledgments}
The author would like to thank Professor Shin Satoh for his helpful comments. The author is supported by JSPS Research Fellowships for Young Scientists (JSPS KAKENHI Grant Number 24$\cdot$9014).

\end{document}